\documentclass[mathpazo]{cicp}
\usepackage{amsmath}
\usepackage{mathrsfs}
\usepackage{graphicx}
\usepackage{subfigure}
\usepackage{rotating}
\usepackage{amssymb}
\usepackage{epsfig}
\usepackage{color}
\usepackage{amsmath}
\usepackage{wrapfig}
\usepackage{threeparttable}
\usepackage{dcolumn}
\newcolumntype{d}{D{.}{.}{-1}}
\usepackage{nomencl}
\makeglossary
\usepackage{fancyvrb}
\fvset{fontsize=\footnotesize,xleftmargin=2em}
\usepackage{lettrine}

\numberwithin{equation}{section}

\begin{document}


\title{Numerical Solution of 3D Poisson-Nernst-Planck Equations Coupled with Classical Density Functional Theory for Modeling Ion and Electron Transport in a Confined Environment}
\author{Da Meng, Bin Zheng, Guang Lin, and Maria L. Sushko\footnote{The first two authors contribute equally to this work. Corresponding authors: G. Lin, Guang.Lin@pnnl.gov, M.L. Sushko, Maria.Sushko@pnnl.gov}}

\address{Pacific Northwest National Laboratory, Richland, WA 99352 USA}



\begin{abstract}

We have developed efficient numerical algorithms for solving 3D steady-state Poisson-Nernst-Planck (PNP) equations with excess chemical potentials described by the classical density functional theory (cDFT). The coupled PNP equations are discretized by a finite difference scheme and solved iteratively using the Gummel method with relaxation. The Nernst-Planck equations are transformed into Laplace equations through the Slotboom transformation. Then, the algebraic multigrid method is applied to efficiently solve the Poisson equation and the transformed Nernst-Planck equations. A novel strategy for calculating excess chemical potentials through fast Fourier transforms is proposed, which reduces computational complexity from $O(N^2)$ to $O(N\log N)$, where $N$ is the number of grid points. Integrals involving the Dirac delta function are evaluated directly by coordinate transformation, which yields more accurate results compared to applying numerical quadrature to an approximated delta function. Numerical results for ion and electron transport in solid electrolyte for lithium-ion (Li-ion) batteries are shown to be in good agreement with the experimental data and the results from previous studies. 


\end{abstract}



\maketitle

\section{Introduction}
\label{sec:Introduction}

Poisson-Nernst-Planck (PNP) equations are widely used to describe the macroscopic properties of ion transport in electrochemical systems \cite{Richardson:2007kx,Bazant:2009vn,Soestbergen:2010ys,Ciucci:2011fk,Marcicki:2012uq} (e.g., lithium-ion (Li-ion) batteries, fuel cells) and biological membrane channels \cite{Eisenberg:1998bh,Kurnikova_BJ1999,Cardenas:2000dq,Hollerbach:2000cr,Coalson:2005nx,Lu:2007vn,Bolintineanu:2009kx,Singer:2009oq}. PNP equations are also known as the drift-diffusion equations for the description of currents in semiconductor devices \cite{Selberherr:1984tg,Markowich:1986ij,Rouston:1990hc,Newman:1991kl}. In these models, excess chemical potential of mobile ions drives their diffusion. However, a highly simplified description of the interactions limited to Coulomb interactions between all charged species, is often used.  To overcome this oversimplification in the representation of collective interactions, classical density functional theory (cDFT) can be used. 
cDFT is a powerful analytical tool to describe mesoscopic interactions, such as excluded volume effects and electrostatic correlation interactions, and thermodynamic properties of inhomogeneous systems from first principles \cite{Singh:1991kx}. 
The PNP-cDFT model is a generalization of the PNP model often used to describe fluids of charged hard spheres in a confined environment. It has been applied to study the selectivity and ionic flux in biological ion channels \cite{Gillespie:2002ly,Gillespie:2005vn,Gillespie:2008ys,Gillespie:2008zr} and shown to provide computational results in good agreement with experimental data and/or theoretical analysis. 


In solid state ion and electron diffusion is also affected by the barriers for elementary transport processes: ion hopping between the adjacent equilibrium sites and electron hopping between the cations in the lattice. Similarly, in biological ion channels, short-range dispersion interactions between the ions and functional groups in the channel proteins would also affect their diffusion. These short-range interactions have a quantum mechanical nature, which makes it challenging to evaluate them analytically. To include these short-range interactions in cDFT model, quantum mechanical simulations can be used to evaluate the barriers for the elementary transport processes and represent the interactions with a square-well potential, featuring depth equal to the barrier and the width comparable to ionic diameters \cite{Sushko_CPL2010,Cao:2005kx}.

To summarize, in our approach --- apart from Coulomb interactions --- electrostatic correlation
and excluded volume effects are treated using cDFT with short-range interactions quantum mechanically evaluated. This approach is equally applicable to study ion and electron transport in nanostructured materials, ion transport through biological ion channels, and small molecule diffusion in mesoporous materials.

In this work, we use the solid electrolyte, lithium phosphorus oxynitride (LiPON), for Li-ion batteries as a test system and study temperature dependence of $\text{Li}^+$ conductivity in LiPON films. This material has a complex $\text{Li}^+$ diffusion pathway \cite{Du:2007fk,Sushko_CPL2011}, which requires a full 3D model for ion and electron transport.  Our previous work showed the PNP-cDFT model's unique capability to capture the physics of nanostructured electrode materials for Li-ion batteries, providing insights into the origin of size effects of conductivity and temperature dependence. The model can be used to guide synthesis of new nanocomposite materials with significantly improved electrochemical properties \cite{Sushko_CPL2010,Sushko_JPCC2010,Sushko_JPCB2010,Sushko_CPL2011,Li:2012uq,Hu:2013fk}.  This research also revealed limitations in modeling realistic nanocomposites with complex structures, calling for optimization of the efficiency of the PNP-cDFT solvers.

The PNP-cDFT model is governed by nonlinear integro-differential equations. The mathematical analysis and numerical simulation generate interesting and challenging problems \cite{Sushko_CPL2011,Gillespie:2002ly,Gillespie:2005vn,Gillespie:2008ys,Gillespie:2008zr,Golovnev:2011qf,Ji:2012kx,Liu:2012ve}.
Numerical methods for solving the PNP system of equations have been studied extensively, including the finite difference \cite{Kurnikova_BJ1999,Cardenas:2000dq}, finite volume \cite{Wu:2002fk,Mathur:2009bs,Bolintineanu:2009kx}, and finite element methods \cite{Hollerbach:2000cr,Lu:2007vn,Coco:2007fv,Lu:2010dz}. In \cite{Zheng:2011fu}, a second-order convergent numerical method is constructed to handle discontinuous dielectric constants and singular sources in the context of biological ion channel applications. 

In the current study, a standard finite difference scheme is sufficient because the diffusion coefficients for PNP equations are assumed to be constant, and, for our applied interests, the computational domain is regular. The main computational challenges for the PNP-cDFT simulations include the solution of large sparse linear systems resulting from the discretizations of PNP equations, and a significant amount of 3D integrals to be calculated for cDFT. In particular, the computation of chemical potentials of charged species requires $O(N^2)$ operations, where $N$ is the number of computational grid points. This makes application of the model to realistic systems computationally very expensive. Hence, most of the existing studies are restricted to the 1D case \cite{Gillespie:2002ly,Gillespie:2005vn,Gillespie:2008ys,Gillespie:2008zr,Golovnev:2011qf,Hyon:2011uq,Ji:2012kx,Liu:2012ve}. Numerical simulations of a 3D PNP-cDFT model have been reported in \cite{Kurnikova_BJ1999,Sushko_CPL2011}, using coarse grid to reduce computational complexity. In this work, we apply a state-of-the-art fast Poisson solver --- algebraic multigrid (AMG) method that has computational complexity of $O(N\log N)$. To speed up numerical integrations, we reformulate them as convolution sums and then employ the fast Fourier transform (FFT) method to reduce the computational complexity to $O(N\log N)$. Some integrals in cDFT calculations involve Dirac delta function. A usual approach is to approximate the Dirac delta function by a smooth function (e.g., Gaussian) and then apply a standard quadrature rule. The additional error introduced by the approximation of the delta function is one drawback to this approach. Here, we use the definition of Dirac delta function and change of variables to transform these 3D integrals into 2D integrals on spheres and remove the singularity in the integrands. Finally, 3D large scale numerical simulations of the PNP-cDFT system are made feasible using the packages BoomerAMG \cite{Falgout:2002fk,Henson_ANM2002} and F3DFFT \cite{Pekurovsky:2012uq}.

The rest of the paper is organized as follows: the multiscale model for nanostructured material and its description through PNP equations and cDFT are presented in Section \ref{sec:gov-eq}. Section \ref{sec:Numeric} describes the numerical methods used in this study, including the finite difference discretization, Gummel iteration with relaxation, AMG solver for sparse linear systems of equations, FFT for calculating excess chemical potential, and the special treatment of integrals involving delta function. The validity, accuracy, and computational complexity of the proposed numerical algorithms are demonstrated in Section \ref{Sec:numerics} via a real application, LiPON film simulation. Finally, Section \ref{sec:Conclusions} presents some concluding remarks.

%


\section{PNP-cDFT model}
\label{sec:gov-eq}

In this section, we offer a detailed description of the PNP model coupled with the cDFT, or PNP-cDFT, and its application to Li-ion batteries simulation.

PNP equations provide a mean-field continuum model for the flows of charged particles in terms of the average density distributions $\rho_i$ and the electrostatic potential $\phi$. The chemical potentials of the charged particles are evaluated by cDFT, which models discrete ion interactions and accounts for particle size effects. In PNP-cDFT theory, the above two models are combined to describe the flow of interacting ion species driven by the excess chemical potentials $\mu_i^{\text{ex}}$.

Li ions diffuse in solids either through hopping between the interstitial sites or vacancy migration mechanisms. There are certain barriers for ion and electron diffusion between equilibrium interstitial sites, which are calculated using quantum mechanical approaches for the corresponding bulk materials \cite{Du:2007fk}. The presence of barriers for $\text{Li}^+ / \text{e}^{-}$ diffusion can be represented by the attraction potential between $\text{Li}^+$ ions or electrons and the corresponding equilibrium sites. The simplest form for such potential is a square-well potential with the well depth equal to the barrier for $\text{Li}^+ / \text{e}^{-}$ hopping between these sites.

\begin{figure}[htbp]
\center
\includegraphics[width=2.5in]{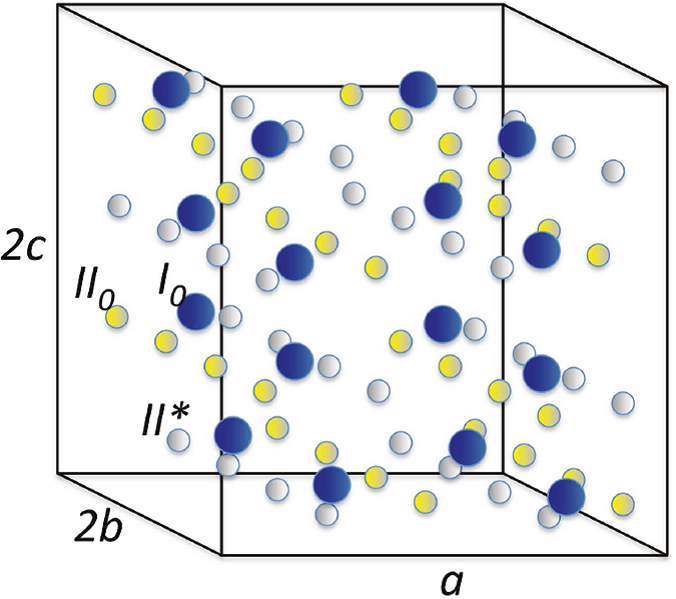}
\caption{Diffusion channel of Li$^+$ in LiPON. Blue spheres are interstitial equilibriums ($I_0$). Yellow spheres ($II_0$) and gray spheres ($II^*$) are metastable sites. a, 2b, and 2c are the sizes along three crystallographic directions.}
\label{fig:ionchannel_model}
\end{figure}

In particular, we focus on the LiPON model, one of the most widely used solid-state electrolytes for thin film batteries developed at Oak Ridge National Laboratory. Figure \ref{fig:ionchannel_model} illustrates a 3D model for the description of $\text{Li}^+$ transport in LiPON. The $I_0$ sites are the equilibrium interstitial sites for $\text{Li}^+$ diffusion (Figure \ref{fig:ionchannel_model}). However, direct $\text{Li}^+$ hopping between $I_0$ sites is energetically unfavorable due to a relatively high energy barrier ($0.21$ eV) and the large distance between these sites ($0.41$ nm), which decrease the probability of ion hopping. According to quantum mechanical DFT simulations, the most energetically favorable path for $\text{Li}^+$ diffusion is a zigzag path via the $II_0$ and $II^*$ sites \cite{Du:2007fk}. Therefore, in our model --- except for site $I_0$ --- intermediate sites $II_0$ and $II^*$ also are introduced (Figure \ref{fig:ionchannel_model}). In our cDFT calculations, these lattice sites $I_0, II_0, II^*$ are modeled by stationary particles.

%
%
%
%
%
%
%

\newpage

\subsection{The 3D steady-state Poisson-Nernst-Planck equations}
\label{}
PNP theory is a continuum electrodiffusion model that represents ion fluxes in terms of density distribution of ion species and potential gradients. It has been widely used in modeling ion transport in biological ion channels or nanocomposites with broad applications in biology and material sciences \cite{Bolintineanu:2009kx,Lu:2007vn,Sushko_CPL2010,Sushko_CPL2011}. 

Within the PNP formalism, the ion flux $\boldsymbol{J}_i$ (for the ions of type $i$, i.e., $\text{Li}^+$ or electrons) in the stationary condition can be calculated in terms of density (or concentration) gradient and potential gradient as follows:
\begin{align}
 &-\boldsymbol{J}_i   = D_i(\boldsymbol{r})\left[
\nabla \rho_i +
\frac{1}{k_BT} \rho_i
\left(
q_i e  \nabla\phi + \nabla\mu^{\text{id}}_{i}(\boldsymbol{r}) +
\nabla\mu^{\text{ex}}_i(\boldsymbol{r})
\right)
\right],\label{NP-1}\\
&\nabla \cdot \boldsymbol{J}_i    = 0, \label{NP-2}\\
&-\nabla\cdot (\epsilon(\boldsymbol{r})\nabla \phi)
= 4\pi\left(\rho_f(\boldsymbol{r}) + e \sum_i q_i \rho_i
\right),
\label{Poisson}
\end{align}
where $\boldsymbol{r} = (x, y, z)$ is the location at which point the functions are defined, $D_i(\boldsymbol{r})$ is the diffusion coefficient of the $i$-th ion species, $\rho_i(\boldsymbol{r})$ is the particle number density, $k_B$ is the Boltzmann's constant, $T$ is the absolute temperature ($k_B T$ is thermal energy), $e$ is the elementary charge, $q_i$ is the valence of species $i$ (with sign), $\phi(\boldsymbol{r})$ is the electrostatic potential, $\epsilon(\boldsymbol{r})$ is the electric permittivity (or dielectric function), $\mu^{\text{ex}}_i(\boldsymbol{r})$ is the excess chemical potential of charged particles at position $\boldsymbol{r}$ (determined within the cDFT framework by (\ref{potential-energy-relation})), $\mu^{\text{id}}_i(\boldsymbol{r})$ is the ideal chemical potential defined in (\ref{c_i_approx}), and $\rho_f(\boldsymbol{r})$ is the fixed charge density in the system. Eqs. (\ref{NP-1}, \ref{NP-2}) may be written in the following form:
\begin{equation}
-\nabla\cdot\left[
D_i(\boldsymbol{r})\left(
\nabla \rho_i+
\frac{1}{k_B T} \rho_i
\left(
q_i e  \nabla\phi + \nabla\mu^{\text{id}}_i(\boldsymbol{r})+
\nabla\mu^{\text{ex}}_i(\boldsymbol{r})
\right)
\right)
\right] = 0.
\label{NP-unsimplified}
\end{equation}
Eq. (\ref{NP-unsimplified}) often is referred to as the drift-diffusion or electrodiffusion equation.


By introducing the effective densities (also referred to as Slotboom variables in semiconductor literature)
\begin{equation}
\bar{\rho}_i =  \rho_i\, e^{(q_i e \phi + \mu_i^{\text{id}} + \mu_i^{\text{ex}})/k_B T},
\label{slotboom-transformation}
\end{equation}
the set of steady-state Nernst-Planck equations (one for each species) can be simplified as:
\begin{equation}
-\nabla\cdot \left(\bar{D}_i \nabla \bar{\rho}_i\right) = 0,
\label{NP-simplified}
\end{equation}
where 
$$
\bar{D}_i  = D_i\, e^{-(q_i e \phi + \mu_i^{\text{id}} + \mu_i^{\text{ex}})/k_B T}.
$$
The Slotboom transformation (\ref{slotboom-transformation}) removes the convection term in Eq. (\ref{NP-unsimplified}) and results in self-adjoint Laplace equations which can be solved efficiently by multigrid method. However, possible large variation of the transformed diffusion coefficient $\bar{D}_i$ could result in large condition number of the stiffness matrix \cite{Lu:2010dz}. Furthermore, the Slotboom transformation may cause overflow problem in its numerical implementation due to the exponential term. Fortunately, this overflow problem does not happen in the applications discussed in this paper.


\subsection{Boundary conditions}
Usually, the PNP equations are accompanied by Dirichlet- and/or Neumann-type boundary conditions \cite{Burger:2012uq}. 

Let the computational domain be given by $\Omega=(0, L_x)\times (0, L_y) \times (0, L_z)$. The external electrostatic  potential $\phi$ is influenced by applied potential, which can be modeled by prescribing Dirichlet boundary condition in $y$-direction as:
\begin{equation}
\phi(\boldsymbol{r}) = \phi_0(\boldsymbol{r}),\;\boldsymbol{r}\in \Gamma_D\subset \partial \Omega,
\end{equation}
where $\Gamma_D=\{(x,y,z)\in\Omega \;|\;y=0 ~\;\text{or}~\; y=L_y\}$. For the remaining part of the boundary $\partial\Omega\backslash \Gamma_D$ (i.e., in $x$ and $z$ directions), a no-flux boundary condition is applied:
\begin{equation*}
\epsilon(\boldsymbol{r})\nabla \phi(\boldsymbol{r})\cdot \boldsymbol{n} = 0,\;\boldsymbol{r}\in \partial \Omega\backslash \Gamma_D.
\end{equation*}
The same types of boundary conditions are imposed for variables $\bar{\rho}_i$ in the transformed Nernst-Planck equations, i.e.,
\begin{equation}
\bar{\rho}(\boldsymbol{r}) = \bar{\rho}_0(\boldsymbol{r}),\;\boldsymbol{r}\in \Gamma_D,
\end{equation}
and
\begin{equation*}
\boldsymbol{J}_i(\boldsymbol{r})\cdot \boldsymbol{n} = 0,\;\boldsymbol{r}\in \partial \Omega\backslash \Gamma_D.
\end{equation*}

The existence and uniqueness of the solution for the nonlinear PNP boundary value problems have been studied in \cite{Park:1997fk,Liu:2005uq,Ji:2012kx} for the 1D case and in \cite{Jerome:1985uq,Burger:2012uq} for multidimensions.

\subsection{Classical density functional theory}
\label{cdft}
cDFT is an analytical tool to evaluate the chemical potentials of charged species modeled by hard-sphere mixtures \cite{Rosenfeld:1989fk,Rosenfeld:1993uq,Cao:2005kx}.  For a given cDFT, an analytical expression for the Helmholtz free energy $F$ is formulated as a functional of the set of particle density distributions for all ion species $\{\rho_i(\boldsymbol{r})\}$. The Helmholtz free energy separates naturally into two terms, the ideal-gas term ($F^{\text{id}}$) that is obtained from classical statistical mechanics
\begin{equation}
F^{\text{id}}[\{\rho_i\}] = k_B T \sum_{i}\int_\Omega \rho_i(\boldsymbol{r})
\left(
\ln (\rho_i(\boldsymbol{r})\Lambda_i^3) - 1 
\right)\,\mathrm{d} \boldsymbol{r},
\end{equation}
where $\Lambda_i$ is the thermal wavelength of component $i$, and the excess free energy ($F^{\text{ex}}$), which has contributions from the internal interactions in the system. Accordingly, the chemical potential can be expressed as the sum of $\mu^{\text{id}}$ and $\mu^{\text{ex}}$. In general, the excess chemical potential is a non-local functional of the ion densities.

The ideal chemical potential is expressed as \cite{Hubbard:2002fk,Wang:2004fk}:
\begin{equation}
\mu_i^{\text{id}}(\boldsymbol{r}) = -\ln \left[  \gamma_i\rho_i(\boldsymbol{r})/\rho_i^{\text{bulk}} \right],
\label{c_i_approx}
\end{equation}
where the activity coefficient $\gamma_i$ is described by the extended Debye-H\"{u}ckel theory \cite{Butler:1998fk,Merkel:2008kx}:
\begin{equation}
\ln \gamma_i = -A q_i^2 \frac{\sqrt{I}}{1 + B a\sqrt{I}}.
\label{activity-coefficient}
\end{equation}
In the preceding formula, $I=\frac{1}{2}\sum_i \rho_i q_i^2$ is the ionic strength, $A=1.82\times 10^6 (\varepsilon T)^{-3/2}$ ($\varepsilon$ is the dielectric constant), and $B=50.3(\varepsilon T)^{-1/2}$.

For a system of charged hard spheres, the excess free energy usually includes contributions from the free energies of Coulomb interactions $F_{\text{C}}^{\text{ex}}$, electrostatic correlations $F_{\text{el}}^{\text{ex}}$, and hard-sphere repulsion $F_{\text{hs}}^{\text{ex}}$ \cite{Patra:1999uq,Li:2004fk}. In \cite{Sushko_JPCB2010}, an additional term corresponding to short-range attraction interactions $F_{\text{sh}}^{\text{ex}}$ was included, resulting in the following decomposition:
\begin{equation}
F^{\text{ex}} = F_{\text{C}}^{\text{ex}} + F_{\text{el}}^{\text{ex}} + F_{\text{hs}}^{\text{ex}} + F_{\text{sh}}^{\text{ex}}.
\end{equation} 

The excess chemical potential can be calculated from the functional derivative of the excess free energy with respect to particle density:
\begin{equation}
\mu_i^{\text{ex}} (\boldsymbol{r})= 
\frac{\delta F^{\text{ex}}[\{\rho_k\}]}{\delta \rho_i}=
\frac{\delta F^{\text{ex}}(\boldsymbol{\rho}(\boldsymbol{r}))}{\delta \rho_i(\boldsymbol{r})},
\label{potential-energy-relation}
\end{equation}
where $\boldsymbol{\rho}=[\rho_1, \dots, \rho_s]$, and $s$ is the number of ion types.

\subsubsection{Hard-sphere component}
The hard-sphere model is often used in statistical mechanics to represent the short-range repulsion between two particles, known as the excluded-volume effect. The fundamental measure theory (FMT) \cite{Rosenfeld:1989fk} and modified fundamental measure theory (MFMT) \cite{Roth:2002fk,Yu:2002fk} are among the most accurate formulations for the description of the structure and thermodynamic properties of inhomogeneous hard-sphere fluids. In this model, the excess Helmholtz free energy functional due to the hard-core repulsion $F_{\text{hs}}^{\text{ex}}$ can be expressed as an integral of the functional of weighted densities, i.e.,
\begin{equation}
F_{\text{hs}}^{\text{ex}}[\{\rho_i\}] = k_B T \int_\Omega \Phi_{\text{hs}}(\boldsymbol{r})\,\rm{d}\boldsymbol{r},
\label{hs-energy}
\end{equation}
where $\Phi_{\text{hs}}$ is a function of weighted densities $n_\alpha$ and $\boldsymbol{n}_\beta$ given by FMT:
\begin{equation}
\Phi_{\text{hs}}(\boldsymbol{r}) =-n_0 \ln (1-n_3)
+\frac{n_1 n_2}{1-n_3}+\frac{n_2^3}{24\pi(1-n_3)^2}
-\frac{\boldsymbol{n}_1\cdot \boldsymbol{n}_2}{1-n_3}-\frac{n_2\boldsymbol{n}_2\cdot\boldsymbol{n}_2}{8\pi(1-n_3)^2}
\label{hard-sphere-weighted-density}
\end{equation}
or MFMT:
\begin{multline}
\Phi_{\text{hs}}(\boldsymbol{r}) =  -n_0\ln(1-n_3)+\frac{n_1n_2}{1-n_3} + 
\left[\frac{1}{36\pi n_3^2}\ln(1-n_3) + \frac{1}{36\pi n_3(1-n_3)^2}
\right] n_2^3\\
- \frac{\boldsymbol{n}_1\cdot\boldsymbol{n}_2}{1-n_3} - \left[
\frac{1}{12\pi n_3^2}\ln(1-n_3)+\frac{1}{12\pi n_3(1-n_3)^2}
\right]n_2(\boldsymbol{n}_2\cdot\boldsymbol{n}_2).
\label{hard-sphere-weighted-density-modified}
\end{multline}
$n_\alpha$ and $\boldsymbol{n}_\beta$ are the weighted average of the density
distribution functions $\rho_i(\boldsymbol{r})$ and are defined by:
\begin{align*}
n_\alpha (\boldsymbol{r}) & = \sum_i \int_\Omega \rho_i(\boldsymbol{r}')\omega_i^{(\alpha)}(\boldsymbol{r}'-\boldsymbol{r})\,{\rm d}\boldsymbol{r}',\;\; \alpha = 0, 1, 2, 3,\\
\boldsymbol{n}_\beta(\boldsymbol{r}) & = \sum_i \int_\Omega \rho_i(\boldsymbol{r}')\boldsymbol{\omega}_i^{(\beta)}(\boldsymbol{r}'-\boldsymbol{r})\,{\rm d}\boldsymbol{r}',\;\;\beta = 1, 2,
\end{align*}
where the ``weight functions" $\omega_i^{(\alpha)}$ and $\boldsymbol{\omega}_i^{(\beta)}$ characterizing the geometry of particles (hard sphere with radius $R_i$ for ion species $i$) are given by:
\begin{align}
\omega_i^{(3)}(\boldsymbol{r})&=\theta (|\boldsymbol{r}|-R_i)\\
\omega_i^{(2)}(\boldsymbol{r})&=|\nabla\theta(|\boldsymbol{r}|-R_i)|=\delta(|\boldsymbol{r}|-R_i)
\label{w2}\\
\boldsymbol{\omega}_i^{(2)}(\boldsymbol{r}) & =\nabla\theta(|\boldsymbol{r}|-R_i)=\frac{\boldsymbol{r}}{r}\delta(|\boldsymbol{r}|-R_i)
\label{wv2}\\
\omega_i^{(0)}(\boldsymbol{r})&=\omega_i^{(2)}(\boldsymbol{r})/(4\pi R_i^2)
\label{w0}\\
\omega_i^{(1)}(\boldsymbol{r})&=\omega_i^{(2)}(\boldsymbol{r})/(4\pi R_i)
\label{w1}\\
\boldsymbol{\omega}_i^{(1)}(\boldsymbol{r}) & = \boldsymbol{\omega}_i^{(2)}(\boldsymbol{r})/(4\pi R_i).
\label{wv1}
\end{align}
In the preceding formulae,  $\theta$ is the Heaviside step function with $\theta(x)=0$ for $x>0$ and $\theta(x) = 1$ for $x\leq 0$, and $\delta$ denotes the Dirac delta function.

From Eq. (\ref{potential-energy-relation}) and Eq. (\ref{hs-energy}), it follows that the hard-sphere chemical potential is given by:
\begin{equation*}
\mu_i^{\text{hs}}(\boldsymbol{r}) = k_B T 
\left(\sum_{\alpha}\int_\Omega \frac{\partial \Phi_{\text{hs}}}{\partial n_\alpha}(\boldsymbol{r}')\omega_i^{(\alpha)}(\boldsymbol{r}-\boldsymbol{r}')\,\rm{d}\boldsymbol{r}'\right)  + 
k_B T \left(\sum_\beta\int_\Omega\frac{\partial \Phi_{\text{hs}}}{\partial \boldsymbol{n}_\beta}(\boldsymbol{r}')\boldsymbol{\omega}_i^{(\beta)}(\boldsymbol{r}-\boldsymbol{r}')\,\rm{d}\boldsymbol{r}'\right).
\end{equation*}

By taking partial derivatives of Eq. (\ref{hard-sphere-weighted-density}) with respect to $n_0$, $n_1$, $n_2$, $n_3$, $\boldsymbol{n}_{1,x}$, $\boldsymbol{n}_{1,y}$, $\boldsymbol{n}_{1,z}$, $\boldsymbol{n}_{2,x}$, $\boldsymbol{n}_{2,y}$, or $\boldsymbol{n}_{2,z}$, we can get
\begin{equation*}
\frac{\partial \Phi_{\text{hs}}}{n_0}  = -\ln (1-n_3)
\end{equation*}

\begin{equation*}
\frac{\partial \Phi_{\text{hs}}}{n_1}  = \frac{n_2}{1-n_3}
\end{equation*}

\begin{equation*}
\frac{\partial \Phi_{\text{hs}}}{n_2}  = \frac{n_1}{1-n_3} + \frac{n_2^2-\boldsymbol{n}_2\cdot\boldsymbol{n}_2}{8\pi (1-n_3)^2}
\end{equation*}

\begin{equation*}
\frac{\partial \Phi_{\text{hs}}}{n_3}  = \frac{n_0}{1-n_3}+
\frac{n_1n_2-\boldsymbol{n}_1\cdot\boldsymbol{n}_2}{(1-n_3)^2}
+\frac{n_2^3-3n_2\boldsymbol{n}_2\cdot\boldsymbol{n}_2}{12\pi(1-n_3)^3}
\end{equation*}

\begin{equation*}
\frac{\partial \Phi_{\text{hs}}}{\partial \boldsymbol{n}_{1,x}} = \frac{\boldsymbol{n}_{2,x}}{n_3-1}
\end{equation*}

\begin{equation*}
\frac{\partial \Phi_{\text{hs}}}{\partial \boldsymbol{n}_{1,y}} = \frac{\boldsymbol{n}_{2,y}}{n_3-1}
\end{equation*}

\begin{equation*}
\frac{\partial \Phi_{\text{hs}}}{\partial \boldsymbol{n}_{1,z}} = \frac{\boldsymbol{n}_{2,z}}{n_3-1}
\end{equation*}

\begin{equation*}
\frac{\partial \Phi_{\text{hs}}}{\partial \boldsymbol{n}_{2,x}} = \frac{\boldsymbol{n}_{1,x}}{n_3-1}-
\frac{n_2 \boldsymbol{n}_{2,x}}{4\pi (1-n_3)^2}
\end{equation*}

\begin{equation*}
\frac{\partial \Phi_{\text{hs}}}{\partial \boldsymbol{n}_{2,y}} = \frac{\boldsymbol{n}_{1,y}}{n_3-1}-
\frac{n_2 \boldsymbol{n}_{2,y}}{4\pi (1-n_3)^2}
\end{equation*}

\begin{equation*}
\frac{\partial \Phi_{\text{hs}}}{\partial \boldsymbol{n}_{2,z}} = \frac{\boldsymbol{n}_{1,z}}{n_3-1}-
\frac{n_2 \boldsymbol{n}_{2,z}}{4\pi (1-n_3)^2}.
\end{equation*}

\subsubsection{Short-range interactions}

In addition to the hard-core repulsion, the short-range attractive particle-particle interactions may be modeled by the following square-well potential:
 \begin{eqnarray*}
\Phi_{\alpha\beta}(r)=\left\{
\begin{array}{cl}
\infty,\;\;& 0\leq r<\sigma_{\alpha\beta}\\
-\epsilon_{\alpha\beta}, \;\;& \sigma_{\alpha\beta} \leq r \leq \gamma \sigma_{\alpha\beta},\\
0,\;\;&r>\gamma\sigma_{\alpha\beta}
\end{array}
\right.
\end{eqnarray*}
where $r$ is the distance between the centers of the spherical particles, $\sigma_{\alpha\beta}=(\sigma_\alpha+\sigma_\beta)/2$ ($\sigma_\alpha$ is the particle hard-core diameter, and $r<\sigma_{\alpha\beta}$ characterizes a hard-core repulsion), $\gamma \sigma_{\alpha\beta}$ is the square-well width, $(\gamma-1)\sigma_{\alpha\beta}$ indicates the range of attraction, $\epsilon_{\alpha\beta}$ is the well depth (positive value represents attractive interaction, while negative value corresponds to repulsive interaction), and the attractive width $\gamma=1.2$ as in reference \cite{Cao:2005kx}. The corresponding mean-field approximation of the free energy is given by:
\begin{equation}
F^{\text{ex}}_{\text{sh}}
=
\frac{1}{2}\int_\Omega\int_\Omega \,{\rm d}\boldsymbol{r}\,{\rm d}\boldsymbol{r}'\,
\sum_{\alpha, \beta=+,s}
\rho_\alpha(\boldsymbol{r})\rho_\beta(\boldsymbol{r}')\Phi_{\alpha\beta}(|\boldsymbol{r}-\boldsymbol{r}'|),
\label{short-range-energy}
\end{equation}
where ``s" denotes the stationary points corresponding to the lattice sites $I_0, II_0, II^*$ shown in Figure \ref{fig:ionchannel_model}.

It follows from (\ref{short-range-energy}) that the short-range chemical potential is given by:
$$
\mu^{\text{sh}}_\alpha (\boldsymbol{r}) = \frac{1}{2}\int_\Omega \sum_\beta \rho_\beta(\boldsymbol{r}')
\Phi_{\alpha\beta}(|\boldsymbol{r}-\boldsymbol{r}'|)\,{\rm d} \boldsymbol{r}'.
$$

Short-range interactions between mobile species (Li$^+$ ions and electrons) are not present in the current study \cite{Sushko_JPCC2010,Sushko_CPL2011}. The density profiles of the stationary points $\rho_s$ are given based on the structure of electrode materials. The depth of the potential well is reasonably equal to the barrier height for $\text{Li}^+$ hopping between two adjacent stationary points and was set according to the quantum mechanical and molecular dynamics data for the barrier heights \cite{Du:2007fk}.

%


The density distributions of the stationary points are reasonably represented by the sum of normalized Gaussian ansatz placed at each corresponding lattice site $\boldsymbol{R}_k$:
$$
\rho_s(\boldsymbol{r}) = \left(\frac{\alpha}{\pi}\right)^{3/2}\sum_k e^{-\alpha|\boldsymbol{r}-\boldsymbol{R}_k|^2}.
$$

\subsubsection{Coulomb interactions}
The free energy of long-range Coulomb interactions is given by:
\begin{equation}
F_{\text{C}}^{\text{ex}} = \frac{k_B T l_B}{2}\sum_{i,j}\int_\Omega\int_\Omega 
\frac{q_i q_j \rho_i(\boldsymbol{r})\rho_j(\boldsymbol{r}')}{|\boldsymbol{r}-\boldsymbol{r}'|}
{\rm d}\boldsymbol{r}\,{\rm d}\boldsymbol{r}',
\label{coulomb-energy}
\end{equation}
where $q_i$, $q_j$ are the valences of the charged species and the Bjerrum length is defined as $l_B=e^2/(4\pi\varepsilon_0\varepsilon k_B T)$. $\varepsilon_0$ is the vacuum permittivity, $\varepsilon$ is the relative dielectric constant of the media, and the sum is over all ion species $i, j$. From (\ref{coulomb-energy}), we can derive the corresponding Coulomb chemical potential:
$$
\mu_i^{\text{C}}(\boldsymbol{r}) =q_i k_B T l_B \sum_j\int_\Omega \frac{q_j \rho_j(\boldsymbol{r}')}{|\boldsymbol{r}-\boldsymbol{r}'|}\,{\rm d}\boldsymbol{r}'.
$$

\subsubsection{Electrostatic correlations}
In most cDFT methods, the excess Helmholtz energy due to  the electrostatic correlations $F^{\text{ex}}_{\text{el}}$ is given by an analytical expression based on the perturbation of a suitably chosen, position-dependent reference bulk fluid \cite{Boda:2002uq,Gillespie:2003fk}. Often, it is described by a second-order functional Taylor expansion in terms of powers of the density fluctuations $\Delta \rho_i (\boldsymbol{r}) = \rho_i(\boldsymbol{r})-\rho_i^{\text{bulk}}(\boldsymbol{r})$ around a reference system with given bulk density profiles $\{\rho_i^{\text{bulk}}(\boldsymbol{r})\}$ \cite{Rosenfeld:1993uq}:
\begin{multline}
F^{\text{ex}}_{\text{el}}[\{\rho_i\}] =
F^{\text{ex}}_{\text{el}}[\{\rho_i^{\text{bulk}}\}]
-k_B T \sum_{i= +, -}  \int_\Omega 
C_i^{\text{(1), el}}(\boldsymbol{r}) \Delta \rho_i(\boldsymbol{r})\, {\rm d} \boldsymbol{r}\\
-\frac{k_B T}{2}\sum_{i,j=+,-} \int_\Omega\int_\Omega C_{ij}^{\text{(2), el}}(|\boldsymbol{r}-\boldsymbol{r}'|)\Delta \rho_i(\boldsymbol{r}) \Delta \rho_j(\boldsymbol{r}')\, {\rm d}\boldsymbol{r} \,{\rm d} \boldsymbol{r}' +O((\Delta \rho(\boldsymbol{r}))^3),
\end{multline}
where the first- and second-order electrostatic direct correlation functions are defined as:
\begin{equation}
 C_i^{\text{(1), el}}(\boldsymbol{r}) = -\frac{1}{k_B T}\frac{\delta F_{\text{el}}^{\text{ex}}}{\delta \rho_i(\boldsymbol{r})}|_{\rho=\rho^{\text{bulk}}} = -\frac{1}{k_B T}\mu_i^{\text{el}}[\{\rho_j^{\text{bulk}}(\boldsymbol{r})\}]
 \label{def-c_i}
\end{equation}
and
\begin{equation}
 C_{ij}^{\text{(2), el}}(|\boldsymbol{r}-\boldsymbol{r}'|) :=
 -\frac{1}{k_B T} \frac{\delta^2 F_{\text{el}}^{\text{ex}}}{\delta\rho_i(\boldsymbol{r})\delta\rho_j(\boldsymbol{r}')},
\end{equation}
where $\mu^{\text{el}}_i$ is the chemical potential of the mobile ions. 

According to the mean spherical approximation (MSA), the second-order direct correlation function $C_{ij}^{\text{(2), el}}$ \cite{Waisman:1972vn,Mier-y-Teran:1990fk,Yu:2004uq}:
\begin{equation}
 C_{ij}^{\text{(2), el}}(r) \approx
\left\{
\begin{array}{cl}
-\frac{q_iq_je^2}{k_B T\varepsilon}\left[
\frac{2B}{\sigma} - \left(\frac{B}{\sigma}\right)^2 r - \frac{1}{r}
\right],\;\;&r \leq \sigma \\
0,\;\;&r > \sigma 
\end{array}
\right.
\label{delta_c_def_new}
\end{equation}
where $r=|\boldsymbol{r}-\boldsymbol{r}'|$ is the distance between two ions, $\sigma_i$ is the diameter, $\sigma=(\sigma_i+\sigma_j)/2$ is the hard-core interaction distance between charged particles $i$ and $j$, and $B$ is given by:
\begin{equation}
B=[\xi+1-(1+2\xi)^{1/2}]/\xi,
\end{equation}
where
$$
\xi^2 = \kappa^2 \sigma^2=\left[\frac{e^2}{\varepsilon_0\varepsilon k_B T}\sum_i \rho_i^{\text{bulk}} q_i^2\right]\sigma^2
$$
and $\kappa$ denotes the inverse Debye screening length.


The electrostatic correlation component of the chemical potential is then given by:
\begin{equation}
\mu^{\text{el}}_i(\boldsymbol{r}) = 
\mu_i^{\text{el}}[\{\rho_k^{\text{bulk}}(\boldsymbol{r})\}]
-k_B T\sum_j\int_\Omega C_{ij}^{\text{(2), el}}(|\boldsymbol{r}- \boldsymbol{r}'|)(\rho_j(\boldsymbol{r}')-\rho_j^{\text{bulk}}) {\rm d}\boldsymbol{r}'.
\end{equation}

\section{Numerical methods}
\label{sec:Numeric}

\begin{figure}[htbp]
\center
\includegraphics[width=3in]{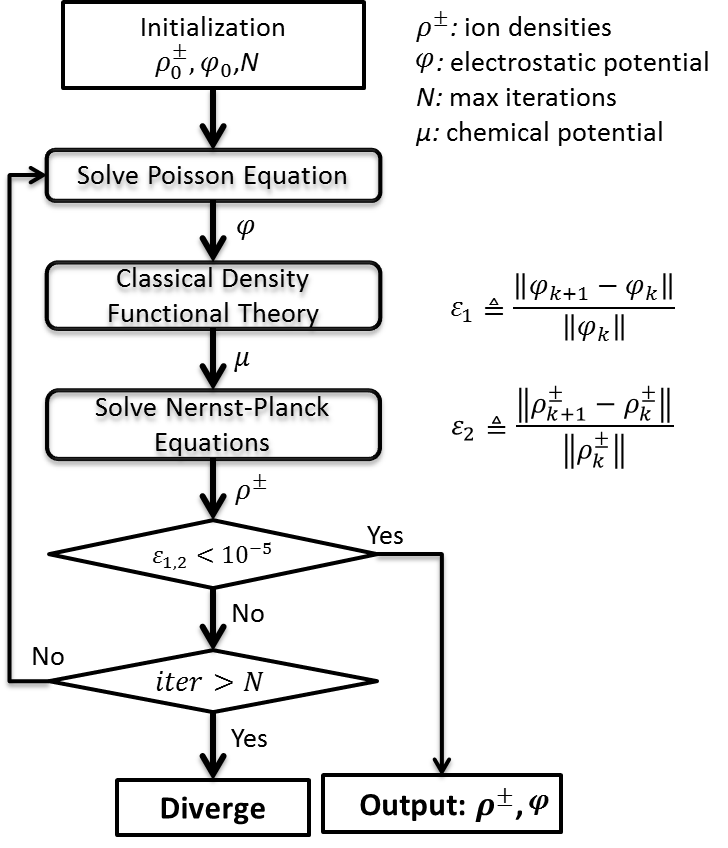}
\caption{PNP-cDFT simulation flow chart.}
\label{fig:pnp-cdft-flow}
\end{figure}

The coupled PNP equations are discretized by finite difference method and solved iteratively using a decoupling method, Gummel iteration \cite{Jerome:1996uq}. The AMG method is applied to solve the Poisson equation (\ref{Poisson}) and the transformed Nernst-Planck equation (\ref{NP-simplified}) efficiently. The excess chemical potential of charged particles are determined by cDFT calculation using FFTs. Figure \ref{fig:pnp-cdft-flow} depicts the flow chart of our numerical simulation. 

We note that the convergence analysis of the numerical method for solving PNP-cDFT system is very difficult and there is no theoretical result to the best of the authors' knowledge. In our numerical convergence test, we are only able to observe that the numerical solutions converge when decreasing the mesh size.
 
\subsection{Finite difference discretization}
The $3$D PNP equations are discretized using the standard $7$-point finite difference scheme on a uniform cubic lattice grid. For example, consider the transformed Nernst-Planck equation (\ref{NP-simplified}) (for simplicity, we drop the species index $i$). The value of $\bar{\rho}$ at a given grid point (corresponding to local index $0$) and its six neighbors (corresponding to the index $j=1, \dots, 6$) satisfy the following equation:
\begin{equation*}
\bar{\rho}^0 \left(\frac{\tilde{D}_{1}+\tilde{D}_{2}}{(\Delta x)^2} +
\frac{\tilde{D}_{3}+\tilde{D}_{4}}{(\Delta y)^2} +
\frac{\tilde{D}_{5}+\tilde{D}_{6}}{(\Delta z)^2}
 \right) = 
 \frac{\bar{\rho}^1\tilde{D}_{1} + \bar{\rho}^2\tilde{D}_{2}}{(\Delta x)^2}
 +\frac{\bar{\rho}^3\tilde{D}_{3} + \bar{\rho}^4\tilde{D}_{4}}{(\Delta y)^2}
 +\frac{\bar{\rho}^5\tilde{D}_{5} + \bar{\rho}^6\tilde{D}_{6}}{(\Delta z)^2}.
\end{equation*}

The diffusion coefficients $\tilde{D}_{j}$ ($j=1, \dots, 6$) are computed via the harmonic mean of the corresponding values at the grid points \cite{Bolintineanu:2009kx}, i.e.,
$$
\tilde{D}_{j} = \frac{2\bar{D}^0\bar{D}^j}{\bar{D}^0+\bar{D}^j},
$$ 
where $\bar{D}^j$ is $\bar{D}$ evaluated at the grid point with local index $j$.

\subsection{Gummel iteration with relaxation}

To obtain the self-consistent solution of the PNP equations, the coupled equations are solved using either the Gummel method with relaxations \cite{Kurnikova_BJ1999} or Newton method \cite{Wu:2002fk,Lu:2007vn}. Gummel iteration is a decoupling method for solving coupled systems of equations. Given an initial guess of the ion concentration profiles $\rho_i$ and the electrostatic potential $\phi$, a new $\phi$ is computed by solving the Poisson equation. Then, the updated $\phi$ is substituted into the Nernst-Planck equation to update the ion concentration. This iterative process terminates when the difference between the results of two subsequent iterations is less than a predefined threshold value ($10^{-6}$ for potential and $10^{-5}$ for ion concentrations in our tests).

The convergence rate of the Gummel iteration is usually slow. To speed up the convergence, successive under- or over-relaxation is employed for the solution updates \cite{Kurnikova_BJ1999,Im_2002fk}. Namely, in the $n$-th PNP iteration, the potential $\phi^{(n)}$ is mixed with the previous potential $\phi^{(n-1)}$ through a relaxation parameter $\lambda_1$ before being substituted into the Nernst-Planck equations,
\begin{equation*}
\phi^{(n)}\leftarrow \lambda_1 \phi^{(n)}+(1-\lambda_1)\phi^{(n-1)}.
\end{equation*}
Similarly, for ion concentrations, we introduce another relaxation parameter $\lambda_2$,
$$
\rho_i^{(n)}\leftarrow \lambda_2 \rho_i^{(n)}+(1-\lambda_2)\rho_i^{(n-1)}.
$$
The relaxation parameters $\lambda_1$ and $\lambda_2$ are selected to achieve rapid convergence while maintaining numerical stability. Because the Nernst-Planck equation is more sensitive to potential change than the Poisson equation to concentration changes, we choose $\lambda_1=0.2$ and $\lambda_2=1.0$ in our numerical simulations. 
  
To solve the coupled nonlinear PNP system, the Newton method, which converges quadratically when a good initial guess is available, has been employed in \cite{Wu:2002fk,Lu:2007vn}.
It requires the construction of a Jacobian matrix that, in practice, may be complicated. In contrast, the Gummel method has a fast initial error reduction, but the convergence rate may be slow. In \cite{Wu:2002fk}, the PNP system is solved with the Newton method, and the resulting linear systems are solved using the generalized minimal residual (GMRES) method with multigrid preconditioning. A better approach may be to combine both methods, i.e., start the solution procedure with a few Gummel iterations to generate a good initial guess then switch to the Newton method to accelerate the convergence.

\subsubsection{The choice of initial guess}
Even when using Gummel method, a good initial guess enhances the stability of solution and speeds up convergence. In our approach, we first solve the static equilibrium problem and determine the distribution of mobile species in external electric potential by minimizing the total free energy with cDFT. These equilibrium density distributions and the corresponding electrostatic potential are subsequently used as initial guesses for the Nernst-Planck and Poisson equations, respectively. This way, we avoid large changes in densities of mobile species when solving the PNP equations by starting from the system's equilibrated state before imposing constant flow conditions. This two-step approach proved to be computationally efficient with typical convergence of the PNP solution in a few tens of iterations.

\subsection{Algebraic multigrid method}
In several early studies of PNP models \cite{Kurnikova_BJ1999,Im_2002fk}, classical Jacobi or Gauss-Seidel iterative methods are used to solve the Poisson and Nernst-Planck equations. These methods suffer from slow convergence, making it difficult to run large-scale 3D simulations. Because PNP equations are all of elliptic type, many existing fast solvers can be applied to solve them efficiently on massively parallel computers. In this work, we use the state-of-the-art fast iterative solver --- AMG method. 

The multigrid (MG) method is well know for its $O(N)$ optimality ($N$ is the number of degrees of freedom) in solving large sparse linear systems resulting from discretizations of partial differential equations. MG's motivation may be described by decomposing the solution error into the sum of high- and low-frequency components. It consists of a smoothing procedure (damped Jacobi, Gauss-Seidel, etc., also called ``smoothers"), which reduces high-frequency error components, and a coarse-grid correction operator for low-frequency errors. Grid transfer operators (restriction and interpolation) are also defined to connect solutions at different grid levels. The success of MG methods relies on the combination of solution procedures at different scales, where different error components are reduced by using the smoothing property of basic iterative methods. For a thorough discussion about the MG method, refer to the books \cite{Hackbusch:1985ys,Wesseling:1992ly,Bramble:1993zr}.

First introduced in the 1980s \cite{Brandt:1982uq,Brandt:1986fk}, the AMG method constructs coarse-grid correction operator and restriction/interpolation operator from a matrix without using grid information. AMG has been quite successful for solving large sparse linear systems, especially those corresponding to discretized elliptic problems. In our numerical simulation, we use BoomerAMG, a parallel implementation of AMG from the hypre (High Performance Preconditioners) library developed at Lawrence Livermore National Laboratory \cite{Falgout:2002fk,Henson_ANM2002}.

\subsection{Fast Fourier transform}
The evaluation of a large number of 3D integrals is another computationally demanding task in cDFT computation (described in Section \ref{cdft}). Direct evaluation of these integrals using numerical quadrature rules is computationally intractable \cite{Burger:2012uq}. Fortunately, these integrals are described in convolution form, which can be calculated effectively using FFT and inverse FFT. In our study, we use P3DFFT, a parallel FFT library optimized for large-scale computer simulations, developed at the University of California, San Diego \cite{Pekurovsky:2012uq}.

More precisely, the 3D integrals involved in the cDFT model can be written in the following convolution form:
\begin{equation}
\int_\Omega  \rho(\boldsymbol{r}') g(|\boldsymbol{r}-\boldsymbol{r}'|)\,{\rm d} \boldsymbol{r}',
\label{3d-convolution-integral}
\end{equation}
where $g(s)=1/s$ for the Coulomb chemical potential, $g(s) = \Phi_{ij}(s)$ for the chemical potential corresponding to the short-range interactions, and $g(s) = C_{ij}^{\text{(2), el}}(s)$ for the electrostatic correlation chemical potential.

In the following, we illustrate how to calculate (\ref{3d-convolution-integral}) using FFT and its inverse on a 1D example. Given a uniform grid on the interval $[0, T]$, $t_i=i\Delta t$, $i=0, 1, \dots, N-1$, where $\Delta t=T/N$. Consider the convolution integral:
\begin{equation}
h(t) = \int_0^T f(\tau) g(|t-\tau|)\,{\rm d} \tau
\label{1d-exmaple-convolution}
\end{equation}
approximated using the trapezoidal rule, i.e.,
\begin{equation}
\begin{array}{rl}
h_k &\approx \sum_{i=0}^{N-1}f_i g_{|k-i|} \Delta x,\;\;\;k=0, 1, \dots, N-1\\
 &= \sum_{i=0}^k f_i g_{k-i} + \sum_{i=k+1}^{N-1} f_i g_{i-k} ,
\end{array}
\label{trapezoidal}
\end{equation}
where $h_k=h(t_k)$, $f_i = f(t_i)$, and $g_{|k-i|}=g(t_{|k-i|})$. Note that for ease of notation, in (\ref{trapezoidal}) we have assumed the integrand has equal values at the two endpoints.

To evaluate Eq. (\ref{trapezoidal}) by applying FFT, we need to write the sums in the convolution form. To this end, we introduce the following auxiliary vectors:
\begin{equation}
\hat{f}_i =\left\{
\begin{array}{cl}
f_i, &\;\;\text{if}\;0\leq i\leq N-1,\\
0, & \;\;\text{if}\;N\leq i\leq 2N-1,
\end{array}
\right.
\end{equation}
\begin{equation}
\hat{g}_i =\left\{
\begin{array}{cl}
g_i, &\;\;\text{if}\;0\leq i\leq N-1,\\
0, & \;\;\text{if}\;N\leq i\leq 2N-1,
\end{array}
\right.
\end{equation}
\begin{equation}
\;\;\;\;\;\;\;\;\;\;\;\;\;\;\tilde{g}_i =\left\{
\begin{array}{cl}
g_{N-2-i}, &\;\;\text{if}\;0\leq i\leq N-2,\\
0, & \;\;\text{if}\;N-1\leq i\leq 2N-1,
\end{array}
\right.
\end{equation}
and define two circular convolution $\hat{h}$ and $\tilde{h}$ (modulo $2N$) as
\begin{align*}
\hat{h}_k &: =(\hat{f}\star\hat{g})_k= \sum_{i=0}^{2N-1} \hat{f}_i \hat{g}_{k-i}, \;\;\forall k\\
\tilde{h}_m &:=(\hat{f}\star \tilde{g})_m = \sum_{i=0}^{2N-1}\hat{f}_i\tilde{g}_{m-i},\;\;\forall m.
\end{align*}
By direct calculations, we get:
\begin{eqnarray*}
&\hat{h}_k =\sum_{i=0}^k f_i g_{k-i},\;\;k=0, 1, \dots, N-1\\
&\tilde{h}_m =\sum_{i=k+1}^{N-1} f_i g_{i-k},\;k=m-(N-1),\;\;m=N-1,\dots, 2N-3,\\
&\tilde{h}_{2N-2} = 0.
\end{eqnarray*}
Hence, 
$$
h_k = \hat{h}_k + \tilde{h}_{k+N-1},\;\;k=0, \dots, N-1,
$$
where $\hat{h}_k$ and $\tilde{h}_m$ can be calculated efficiently by FFT and inverse FFT \cite{Cooley:1967fk} with $O(N\log N)$ computational complexity.

The generalization of the preceding approach for evaluating 3D convolution integrals (\ref{3d-convolution-integral}) is straightforward. As such, the details are omitted here.

\subsection{Integrals involving Dirac delta function}
When calculating hard-sphere chemical potentials, it can be seen from Eqs. (\ref{w2}, \ref{wv2}, \ref{w0}, \ref{w1}, and \ref{wv1}) that several 3D integrals involving Dirac delta functions need to be evaluated. A straightforward approach is to approximate the Dirac delta function by the Gaussian function:
\begin{equation}
\delta(r) \approx \delta_{\alpha}(r) = \frac{1}{\alpha \sqrt{\pi}}e^{-r^2/\alpha^2}
\end{equation}
and then apply any numerical quadrature rule. However, a more efficient and accurate approach is to reduce the 3D volume integrals into 2D spherical integrals using the definition of Dirac delta function and change of variables. More precisely, we have:
\begin{align*}
I(\boldsymbol{r}) &= \int_\Omega \rho(\boldsymbol{r}')\delta(|\boldsymbol{r}'-\boldsymbol{r}|-R)\,{\rm d}\boldsymbol{r}'\\
& = \int_0^{2\pi}\int_0^{\pi}\rho(x+R\cos\theta\sin\varphi, y+R\sin\theta\sin\varphi, z+R\cos\varphi)\\
& \;\;\;\;\;\;\;\;\;\;\;\;\;\;\;\;\;\;\;
R^2 \sin\varphi\,{\rm d}\varphi \,{\rm d}\theta\\
& \approx \sum_{i=1}^k \sum_{j=1}^k
[
\rho(x+R\cos\theta_i\sin\varphi_j, y+R\sin\theta_i\sin\varphi_j, z+R\cos\varphi_j)
\\
& \;\;\;\;\;\;\;\;\;\;\;\;\;\;\;\;\;\;\;
R^2\sin\varphi_j ]\,w_i\,\tilde{w}_j,
\end{align*}
where $\{(\theta_i, \varphi_j)\}$ are the quadrature points and $\{w_i, \tilde{w}_j\}$ are the corresponding weights. Zero extension of the integrand is used when part of the sphere lies out of the computational domain $\Omega$.

%

\section{Numerical results}
\label{Sec:numerics}

We use the Li-ion and electron conducting solid electrolyte, LiPON, to demonstrate our PNP-cDFT solver's performance.

\subsection{Computational domain and physical parameters}

The computational domain is given by $2a\times Mb \times 5c$, where $a=1.053$ nm, $b=0.612$ nm, and $c=0.493$ nm are the lattice parameters for LiPON and $M=10, ..., 200$. $\text{Li}^+$ ions are represented as spherical particles with charge $\text{q}_+=1$ and diameter $\sigma_+=0.06$ nm. Similar representation is used for electrons, diffusing along with $\text{Li}^+$, with the parameters $q_-=-1$ and diameter $\sigma_-=0.001$ nm. The experimental value for the LiPON dielectric constant, 16.6, has been used \cite{Fu:2003fk}. 
The parameters used in our PNP-cDFT simulation for $\text{Li}^+$/electron transport are listed in Table \ref{table:nonlin}.

\begin{table}[ht]
\centering 
\caption{Parameters for computation} 

\begin{tabular}{l l l l} 
\hline 
Parameter & Symbol & Value & Unit \\ [0.5ex] 
\hline 
Diffusion coefficient ($\text{Li}^+$) & $D_+$ & $1 \times 10^{-6}$ & $\text{cm}^2/s$ \\ 
Diffusion coefficient (electron) & $D_-$ & $1\times 10^{-6}$ & $\text{cm}^2/s$ \\
Valence ($\text{Li}^+$) & $q_+$ & $+1$ &  \\
Valence (electron) & $q_-$ & $-1$ &  \\
Dielectric constant of media & $\epsilon$ & $16.6$ &  \\ 
Fixed charge density (at $z=0$) & $\rho_f$ & $0.1$ & $1/\text{nm}^2$ \\
Bulk density of salt & $\rho_{\text{bulk}}$ & $6.02\times 10^{-5}$ & $1/\text{nm}^3$\\
Thermodynamic beta ($T=223K$) & $\beta=\frac{1}{k_B T}$ & $52.03749$ & $1/\text{eV}$\\
LiPON nanoparticle size in $a$-direction & $L_x$ & $2\times 1.053$ & nm\\
LiPON nanoparticle size in $b$-direction & $L_y$ & $150\times 0.612$ & nm\\
LiPON nanoparticle size in $c$-direction & $L_z$ & $5\times 0.493$ & nm\\
Sphere diameter ($\text{Li}^+$) & $\sigma_+$ & $0.06$ & nm\\
Sphere diameter (electron) & $\sigma_-$ & $0.001$ & nm\\
Sphere diameter (stationary points) & $\sigma_s$ & $0.2$ & nm\\
Square-well potential depth ($\text{Li}^+$, $s_1$) & $\epsilon_{+,s_1}$ & $0.21$ & eV\\
Square-well potential depth ($\text{Li}^+$, $s_2$) & $\epsilon_{+,s_2}$ & $0.17$ & eV\\
Square-well potential depth ($\text{Li}^+$, $s_3$) & $\epsilon_{+,s_3}$ & $0.17$ & eV\\
Bjerrum length ($T=223K$) & $l_B = \frac{e^2}{4\pi \epsilon_0\epsilon k_B T}$ & $4.514$ & nm\\
[1ex] 
\hline 
\end{tabular}
\label{table:nonlin} 
\end{table}

\subsection{Results and discussion}

\begin{figure}[htbp]
\center
\includegraphics[width=4in]{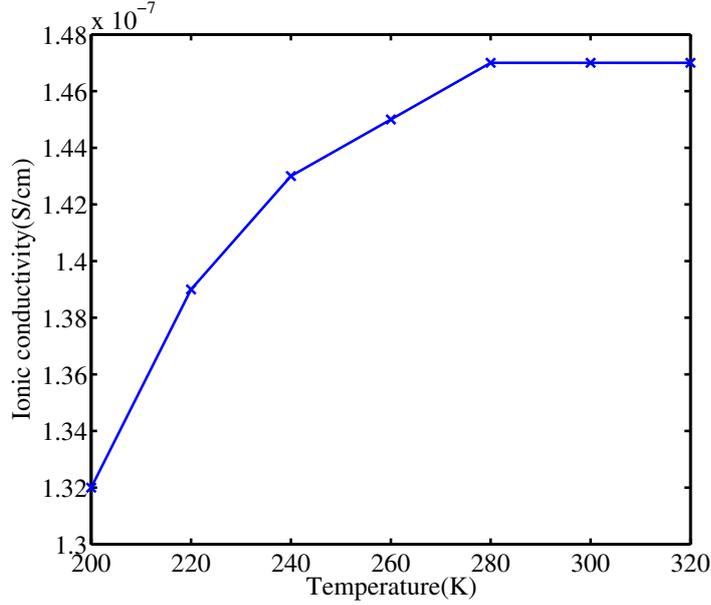}
\caption{$\text{Li}^+$ conductivity in LiPON nanoparticles along the $y$ direction with varying temperature. The particle size in $y$ direction is $150b$ ($91.8$ nm).}
\label{figure:lipon_conductivity_temperature}
\end{figure}

Ion conductivity obtained from numerical simulation is on the order of $10^{-7}$ S/cm in $y$ direction at temperature $298~K$. This is close to the conductivity of $3.0\times 10^{-7}$ S/cm observed experimentally in \cite{Wang:1995vn} for $\text{Li}_{0.99}\text{PO}_{2.55}\text{N}_{0.30}$ glass. Moreover, based on our simulation, conductivity increases from $1.32\times 10^{-7}$ to $1.47 \times 10^{-7}$ when temperature increases from $200~K$ to $320~K$ (Figure \ref{figure:lipon_conductivity_temperature}). This temperature-dependence trend for conductivity agrees well with experimental data \cite{Ribeiro:2012fk}. It is good to note that reproducing the exact experimental value for conductivity is highly improbable as conductivity strongly depends on the film microstructure (shown in \cite{Ribeiro:2012fk}). The simulation domain we use is finite with the size in the main conductivity direction of $61.2 - 91.8$ nm, which corresponds to polycrystalline LiPON films, and the simulated conductivities are well within the limits of variations in experimental values for polycrystalline samples. Apart from film microstructure, the value for the elementary diffusion coefficient, a pre-factor to the conductivity, is another source of uncertainly. In our simulations, we use the value of $10^{-6}\;\text{cm}^2/s$ measured for rutile titanium dioxide, or $\text{TiO}_2$, a material with an interstitial $\text{Li}^+$ diffusion mechanism and similar barriers for elementary diffusion processes \cite{Johnson:1964fk,Koudriachova:2003uq}. Therefore, we expect the rutile value to be a good estimate for the $\text{Li}^+$ diffusion coefficient in LiPON.
%

\begin{figure}[htbp]
\center
\includegraphics[width=4in]{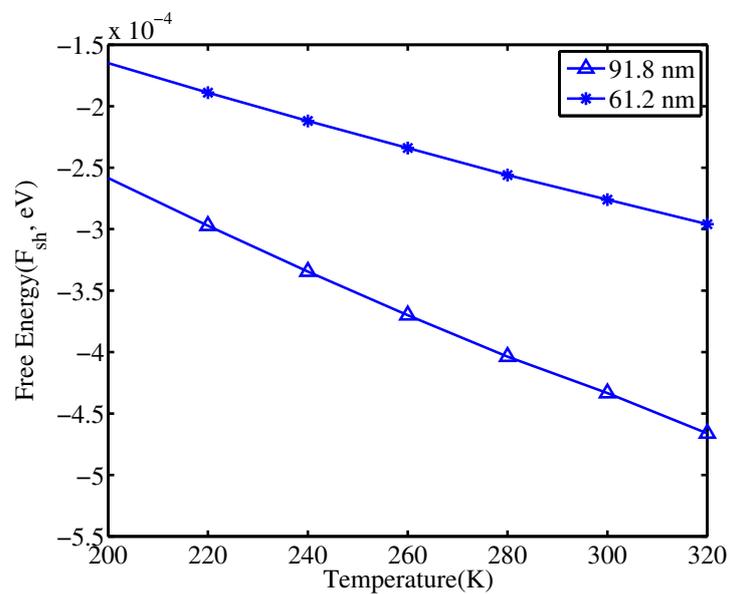}
\caption{Temperature dependence of the $\text{Li}^+$ short-range free energy in LiPON nanoparticles. The particle size in $b$ direction is $100b$ ($61.2$ nm) and $150b$ ($91.8$ nm).}
\label{figure:lipon_shortrange_temperature}
\end{figure}

In our previous work, we have shown that the mechanism of coupled ion and electron transport changes from strongly coupled ionic and electronic to predominantly ionic as the particle size becomes larger than the Debye length. These differences reflect the changes in the mechanism for the compensation of the external electric field, which is achieved through highly correlated ion and electron flux in small nanoparticles and by the formation of the space-charge zone at the surface of large nanoparticles. These different mechanisms for conductivity also lead to different temperature dependences. For large nanoparticles considered here, conductivity increases with temperature due to two reasons: i) thermal motion leads to partial destabilization of the space-charge layer, supplying more ions and electrons to the flux through the nanoparticle, and ii) lowering of the effective barrier for elementary ion transport between adjacent stationary points. The second effect is manifested in the monotonic decrease in the calculated short-range energy $F^{\text{ex}}_{\text{sh}}$, which dominates the total excess free energy $F^{\text{ex}}$ (Figure \ref{figure:lipon_shortrange_temperature}). This effect is weaker for smaller nanoparticles due to the change in relative contributions of short-range and electrostatic correlation free energies.


\begin{figure}[htbp]
\center
\includegraphics[width=4in]{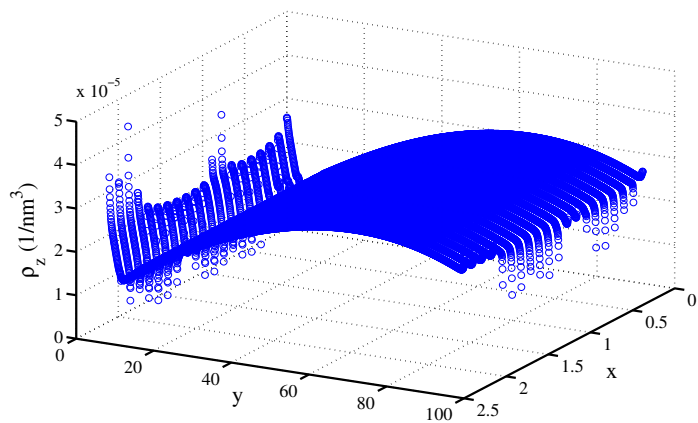}
\caption{$\text{Li}^+$ ion density in $xy$ plane ($z=1.2577$ nm).}
\label{fig:Density-in-z}
\end{figure}

\begin{figure}[htbp]
\center
\includegraphics[width=4in]{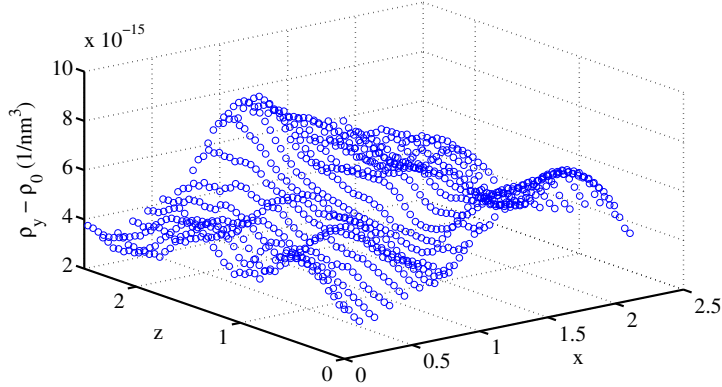}
\caption{$\text{Li}^+$ ion density in $xz$ plane ($\rho_0 = 2.941152817\times 10^{-5}$, $y=55.178$ nm).}
\label{fig:Density-in-y}
\end{figure}

The density distribution in the conduction plane (plane parallel to the $y$  direction) reveals the formation of the space-charge layer at the boundaries of the nanoparticle (Figure \ref{fig:Density-in-z}). In contrast, there is almost no variation in ion and electron densities in the plane normal to the conduction plane with the variations in the density on the order of $10^{-15}$ $\text{nm}^{-3}$ (Figure \ref{fig:Density-in-y}). 

\subsubsection{Size dependency of the $\text{Li}^+$ ion conductivity}

\begin{figure}[htbp]
\center
\includegraphics[width=4in]{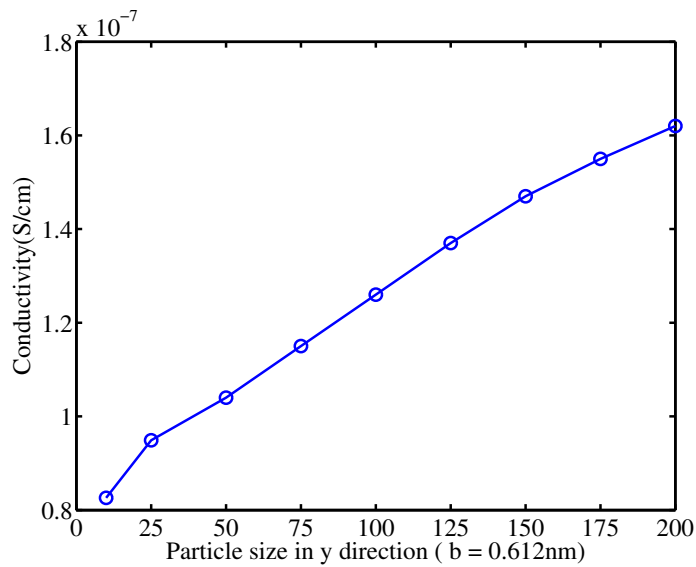}
\caption{Size effect of the conductivity at $300~K$.}
\label{figure:lipon_size_effect}
\end{figure}

The size effects on $\text{Li}^+$ conductivity at temperature $T=300~K$ is shown in Figure \ref{figure:lipon_size_effect}. The observed monotonic increase in conductivity is due to a combination of several competing effects. On one hand, the increase in nanoparticle size leads to the decrease in the effective gradient of electric potential or local electric field in $b$ direction, resulting in stronger correlation of ion and electron fluxes. On the other hand, it reduces the driving force acting on ions and electrons, reducing their diffusivities through the nanoparticle. Overall, the size dependence of the conductivity can be expressed as \cite{Luo:2009fk}:
$$
\sigma_c = \frac{1}{E_y L_y}\int_0^{L_y} J_b(y)\,\mathrm{d} y,
$$
where $\sigma_c$ is the conductivity, $E_y$ is the local electric field, $L_y$ is the nanoparticle's size along the $y$ direction, and $J_b$ is the flux along the $y$ direction. Competition from the previously described effects leads to an almost linear increase in conductivity with nanoparticle sizes ranging from 10 to 100 nm, as observed in our simulations, to very weak dependence of the conductivity on particle sizes ranging from 100 to 1000 nm. Overall, our simulations demonstrate excellent agreement with experimental data and analytic theories, validating our approach \cite{Wang:1995vn,Ribeiro:2012fk}.

\subsubsection{Computational complexity}
\begin{figure}[htbp]
\center
\includegraphics[width=4in]{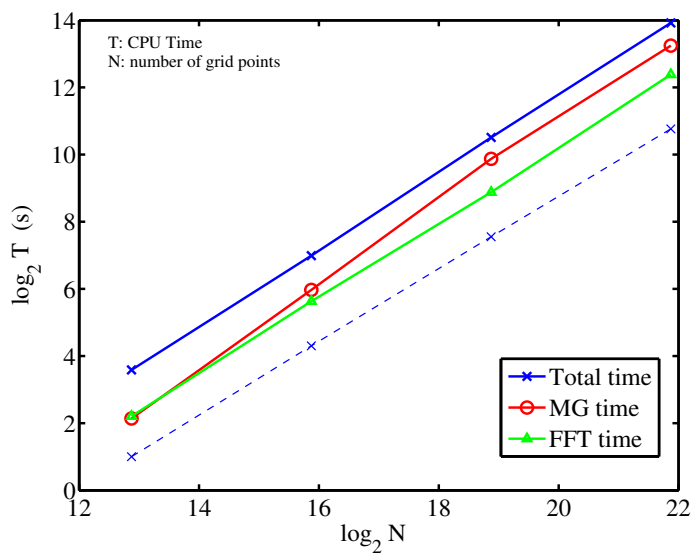}
\caption{Computational complexity (log scale).}
\label{figure:complexity}
\end{figure}

We verify the computational complexity of the methods discussed in Section \ref{sec:Numeric} via numerical experiments. In Figure \ref{figure:complexity}, we plot the graph of the solution time (in seconds), as well as the time for AMG and FFT components versus the number of grid points $N$. Here, we consider a nanoparticle size $2a\times 10b \times 5c$ and run simulations on a sequence of uniformly refined meshes using a single processor. The number of PNP iterations is stable with respect to the mesh size. Hence, the CPU time is a good measurement for the computational complexity of FFT, which is $O(N\log N)$ (see Figure \ref{figure:complexity} where the dashed line represents $O(N\log N)$). However, we do not observe $O(N)$ multigrid complexity. This nonuniform convergence of the multigrid is possibly due to the PNP-cDFT integral-differential system's high nonlinearity.

\section{Conclusions}
\label{sec:Conclusions}
As part of this effort, we have examined numerical methods for simulating a PNP-cDFT model using state-of-the-art AMG and FFT packages. To evaluate chemical potentials efficiently, a novel treatment of 3D integrals via FFT and integrals involving delta function is proposed. The computational complexity of our simulation is $O(N\log N)$, which makes large-scale 3D PNP-cDFT simulation feasible.
Numerical results are validated through comparison with experimental data and results from previous studies. In future work, we will apply the proposed methods to simulate time-dependent PNP-cDFT systems. In addition, parallel scalability of these methods will be reported in a forthcoming paper.
\section*{Acknowledgments}
\label{sec:Acknowledgements}

Work by MLS and DM was supported by the Materials Synthesis and Simulation across Scales (MS3) Initiative (Laboratory Directed Research and Development (LDRD) Program) at Pacific Northwest National Laboratory (PNNL). Work by GL was supported by the U.S. Department of Energy (DOE) Office of Science's Advanced Scientific Computing Research Applied Mathematics program and work by BZ by Early Career Award Initiative (LDRD Program) at PNNL. PNNL is operated by Battelle for the DOE under Contract DE-AC05-76RL01830. The research was performed using PNNL Institutional Computing, as well as the National Energy Research Scientific Computing Center at Lawrence Berkeley National Laboratory.
%
%

\bibliographystyle{elsart-num}
\bibliography{pnp_cdft}

\end{document}